\documentclass[12pt]{article}

\usepackage{amsmath}
\usepackage{amssymb}
\usepackage{amscd}
\usepackage{eucal}
\usepackage{theorem}

\usepackage[dvips]{graphicx}

\title{Hilbert schemes of finite abelian group orbits and Gr\"obner fans}
\author{Tomohito Morita}
\date{}

\makeatletter
\renewcommand\section{%
\@startsection {section}{1}{\z@}%
    {-3.5ex \@plus -1ex \@minus -.2ex}%
    {2.3ex \@plus.2ex\@afterindenttrue}%
    {\normalfont\Large\bfseries}} %
\makeatother

\newcommand{\itfamily}{\itshape}

\theorembodyfont{\rmfamily}
\theoremstyle{change}
\newtheorem{defn}{Definition}[section]

\theorembodyfont{\itfamily}
\theoremstyle{change}
\newtheorem{thm}[defn]{Theorem}

\theorembodyfont{\itfamily}
\theoremstyle{change}
\newtheorem{prop}[defn]{Proposition}

\theorembodyfont{\itfamily}
\theoremstyle{change}

\theorembodyfont{\itfamily}
\theoremstyle{change}
\newtheorem{cor}[defn]{Corollary}

\theorembodyfont{\rmfamily}
\theoremstyle{change}

\theorembodyfont{\rmfamily}
\theoremstyle{change}
\newtheorem{eg}[defn]{Example}

\theorembodyfont{\rmfamily}
\theoremstyle{change}

\theorembodyfont{\rmfamily}
\theoremstyle{change}
\newtheorem{rem}[defn]{Remark}

\theorembodyfont{\rmfamily}
\theoremstyle{change}

\newcommand{\rom}[1]{\textup{#1}}
\let\refsave=\ref
\def\ref#1{\textup{\refsave{#1}}}

\makeatletter
\DeclareRobustCommand{\qed}{%
  \ifmmode 
  \else \leavevmode\unskip\penalty9999 \hbox{}\nobreak\hfill
  \fi
  \quad\hbox{\qedsymbol}}
\newcommand{\openbox}{\leavevmode
  \hbox to.77778em{%
  \hfil\vrule
  \vbox to.675em{\hrule width.6em\vfil\hrule}%
  \vrule\hfil}}
\newcommand{\qedsymbol}{\openbox}
\newenvironment{proof}[1][\proofname]{\par
  \normalfont
  \topsep6\p@\@plus6\p@ \trivlist
  \item[\hskip\labelsep\itshape
    #1\@addpunct{.}]\ignorespaces
}{%
  \qed\endtrivlist
}
\newcommand{\proofname}{Proof}

\DeclareRobustCommand{\fin}{%
  \ifmmode 
  \else \leavevmode\unskip\penalty9999 \hbox{}\nobreak\hfill
  \fi
  \quad\hbox{\finsymbol}}

\newcommand{\finsymbol}{\rule{1.5mm}{3mm}}

\newenvironment{pf}{\proof[Proof]}{\endproof}
\newenvironment{pf*}[1]{\proof[#1]}{\endproof}

\newcommand{\A}{{\mathbb{A}}}
\newcommand{\N}{{\mathbb{N}}}

\newcommand{\Z}{{\mathbb{Z}}}
\newcommand{\R}{{\mathbb{R}}}
\newcommand{\C}{{\mathbb{C}}}
\renewcommand{\P}{{\mathbb{P}}}

\newcommand{\ini}{\operatorname{in}}
\newcommand{\Hilb}{\operatorname{Hilb}}
\newcommand{\Ch}{\operatorname{Ch}}
\newcommand{\str}{\mathcal{O}}

\begin{document}
\maketitle
\begin{center}
Department of Mathematics\\
Tokyo Institute of Technology\\
Oh-okayama, Meguro-ku, Tokyo 152-8551\\
Japan

\

morita.t.ae@m.titech.ac.jp

\end{center}
\begin{abstract}
Let $G$ be a finite abelian subgroup of
 $PGL(r-1,K)=\mathrm{Aut}(\P^{r-1}_K)$. In this paper, we prove that the
 normalization of the $G$-orbit Hilbert scheme $\Hilb^G(\P^{r-1})$ is
 described as a toric variety, which corresponds to the Gr\"obner fan
 for some homogeneous ideal $I$ of $K[x_1, \ldots ,x_r]$.
\end{abstract}

\section{Introduction}
Let $K$ be an algebraically closed field and $G$ a finite subgroup of $GL(n,K) \subset \mathrm{Aut}(\A^n)$ of order prime to the characteristic of $K$. The $G$-orbit Hilbert scheme (or the Hilbert scheme of $G$-orbits) $\Hilb^G(\A^n)$ is introduced by Ito and Nakamura \cite{IN2}. The $G$-orbit Hilbert scheme is the scheme parameterizing all $G$-invariant smoothable zero-dimensional subschemes of $\A^n$ of length $m:= | G | $. Here a smoothable zero-dimensional subscheme $Z$ of $\A^n$ of length $m$ is the subscheme of $\A^n$ for which $H^0(Z,\str_Z)$ is a regular representation of $G$.

Nakamura \cite{Nak} proved that for a finite abelian group $G \subset GL(n,K)$, the normalization of $\Hilb^G(\A^n)$ is the toric variety corresponding to the fan which is defined by the $G$-graph, which is defined in \cite{Nak}, and he also proved that $\Hilb^G(\A^3)$ is a crepant smooth resolution of $\A^3/G$ if $G \subset SL(3,K)$. Furthermore, Ito \cite{Ito} proved the following result. If $G$ is a finite cyclic group of $GL(2,\C)$ and acts freely on $\C^2\backslash\{0\}$, then $\Hilb^G(\C^2)$ is described by the Gr\"obner fan for the ideal $I$ of $\C[x,y]$, corresponding to a subscheme which is a free $G$-orbit contained in $(\C^*)^2$.

In this paper, we consider the case when $G$ is a finite abelian subgroup of the diagonal subgroup of $PGL(r-1,K)=\mathrm{Aut}(\P^{r-1})$. The $G$-orbit Hilbert scheme $\Hilb^G(\P^{r-1})$ defined as before, and the normalization of it is proved to be the toric variety corresponding to the Gr\"obner fan for a homogeneous ideal $I$ of $K[x_1, \ldots ,x_r]$, corresponding to a subscheme which is  a free $G$-orbit contained in $(K^*)^{r-1} \cong \{ x_1 \cdots x_r \neq 0 \} \subset \P^{r-1}$ (Theorem \ref{main}). The corresponding results on $\Hilb^G(\A^{r-1})$ are easily deduced from it (Corollary \ref{HilbA}). This gives an alternative proof and a generalization of Ito's result.

 The proof of our theorem consists of three steps. Note that this diagonal subgroup $T'$ of $PGL(r-1,K)$ is isomorphic to an algebraic torus $(K^*)^{r-1} \subset \P^{r-1}$. Let $\Hilb^m(\P^{r-1})$ denote the Hilbert scheme of $m$ points in $\P^{r-1}$, where $m=|G|$. For any homogeneous ideal $I$ as in the theorem, $I$ defines a point $P$ of $\Hilb^m(\P^{r-1})$, and we prove that the closure in $\Hilb^{|G|}(\P^{r-1})$ of $T'$-orbit of $P$ is coincides with $\Hilb^G(\P^{r-1})$.
Next, we show that the normalization of $\Hilb^G(\P^{r-1})$ is described by the state polytope of $I$, defined in (\ref{state}) (Proposition \ref{Oda}). Finally, we show that the normal fan of the state polytope coincides with the Gr\"obner fan for $I$ by a theorem of Bayer and Morrison \cite{BM} (Theorem \ref{gfan}). 

Gr\"obner fans are computable (for example by \cite{Jen}). Hence we hope that the results in this paper are useful for the study on the $G$-orbit Hilbert schemes, especially in  higher dimensional cases.

After we wrote up this paper, the author found the paper \cite{CMT2} by Craw, Maclagan, and Thomas, where they study the moduli space of McKay quiver representations. Our work is related to their theory in the case of the $G$-orbit Hilbert schemes, but the method is different. T. Yasuda also obtained related results in \cite{Yas}, but his method is different from ours.  Y. Ito communicated to the author that Y. Sekiya independently obtained a similar result to ours.
\bigskip 

\noindent {\bf Acknowledgments.}
I would like to thank Professor Takao Fujita for valuable advice and for pointing out some mistakes in English. I also would like to thank Professors Takeshi Kajiwara and Chikara Nakayama for helpful discussions and comments on earlier versions.

\section{Gr\"obner basis and Gr\"obner fan}

We recall the basic notations for Gr\"obner basis and Gr\"obner fan. These notations are based on the book by Sturmfels \cite{Stu}.\\

Let $K$ be any field and $K[X] =K[x_{1},\ldots,x_{n}]$ the polynomial ring in $n$ indeterminates.
Let $\N$ be the set of non-negative integers. For $a=(a_1,\ldots,a_n) \in \N^n$, let $X^a$ denote the monomial $x_1^{a_1} \cdots x_n^{a_n}$. By this correspondence $a \mapsto X^a$, the lattice $\N^n$ is embedded in $K[X]$ as a multiplicative semigroup and this image coincides with the set of monomials in $K[X]$.

\begin{defn}
A total order $<$ on $\N^n$ is called a \textit{term order} if $(0,\ldots,0)$ is the unique minimal element, and $a<b$ implies $a+c<b+c$ for all $a,b,c \in \N^n$.
\end{defn}

We denote by $\R_+$ the set of non-negative real numbers. 
Fix a weight vector $w=(w_1,\ldots  ,w_n)\in \R^n$ and a term order $<$. Let $\left<\ ,\ \right>$ denote the standard inner product of $\R^n$.

\begin{defn}
For a polynomial $f=\Sigma_{c_i \neq 0} c_iX^{a_i}$, the \textit{initial term} $\ini_w (f)$ of $f$ with respect to $w$ is the sum of all terms $c_iX^{a_i}$ such that the inner product $\left<w,a_i\right>$ is maximal.
For an ideal $I$ of $K[X]$, the \textit{initial ideal} $\ini_{w}(I)$ is the ideal generated by the initial terms of all elements $f$ in $I$:
$$\ini_w(I)=\left<\{\ini_w(f)|f \in I\}\right>.$$
\end{defn}

This ideal need not be a monomial ideal. However, $\ini_w(I)$ is a monomial ideal if $w$ is chosen sufficiently generic for a given $I$.

\begin{defn}
Let $w$ be a weight vector in $\R_+^n$ and $<$ the fixed term order. The \textit{weight term order} $<_w$ is the term order defined as follows:
$$X^a <_w X^b \Leftrightarrow \left<a,w\right> < \left<b,w\right>,\ \mathrm{or}\  \left<a,w\right>=\left<b,w\right> \ \mathrm{and} \  a < b.$$
\end{defn}

\begin{defn}
For a polynomial $f$, the \textit{initial term} $\ini_{<_w}(f)$ of $f$ with respect to $w$ is the maximal term of $f$ with respect to $<_w$.
For an ideal $I$ of $K[X]$, the \textit{initial ideal} $\ini_{<_w}(I)$ is the ideal generated by the initial terms of all elements $f$ in $I$:
$$\ini_{<_w}(I)=\left<\{\ini_{<_w}(f)|f \in I\}\right>.$$
\end{defn}

\begin{defn}
A finite subset $\mathcal{G}\subset I$ is a \textit{Gr\"obner basis} if $\ini_{<_w}(I)$ is generated by $\{\ini_{<_w}(g)|g\in \mathcal{G}\}$.
\end{defn}

It is known that $I$ is generated by the Gr\"obner basis of $I$.

\begin{defn}
The Gr\"obner basis $\mathcal{G}$ of $I$ is \textit{reduced} if for any two elements $g,h\in \mathcal{G}$, no term of $h$ is divisible by $\ini_{<_w}(g)$. 
\end{defn}

It is known that the reduced Gr\"obner basis of $I$ is unique.

\begin{defn}
Two weight vectors $w,w' \in \R^n$ are called \textit{$I$-equivalent} (or simply \textit{equivalent}) if $\ini_w(I)=\ini_{w'}(I)$.
\end{defn}

Then we can consider the equivalence classes of weight vectors.

\begin{prop}[\textrm{\cite[Proposition 2.3]{Stu}}]
For any ideal $I$, and for any weight vector $w$, the equivalence class $c[w]$ of weight vectors is relatively open convex polyhedral cone. Moreover, if $w$ is contained in $\R_+^n$ and chosen sufficiently generic, $c[w]$ is given by the reduced Gr\"obner basis $\mathcal{G}$ of $I$ with respect to $<_w$ as follows\rom{:}
$$c[w]=\{w' \in \R^n | \ini_{w'}(g)=\ini_w(g) \mathrm{\ for\ any\ } g \in \mathcal{G}\}.$$
\end{prop}

\begin{defn}
The \textit{Gr\"obner region} $GR(I)$ for $I$ is the set of all $w \in \R^n$ such that $\ini_w(I)=\ini_{w'}(I)$ for some $w' \in \R_+^n$. Clearly, $GR(I)$ contains $\R_+^n$.
\end{defn}

\begin{rem}
In general, the $GR(I)$ does not coincide with $\R^n$. However, when the ideal $I$ is homogeneous, it is known that $GR(I)$ coincides with $\R^n$.
\end{rem}

\begin{defn}
The \textit{Gr\"obner fan} $GF(I)$ for an ideal $I$ is the set consisting of the faces of the closed cones of the form $\overline{c[w]}$ for  some $w\in \R^n_+$.
A closed cone $\sigma$ is called \textit{Gr\"obner cone} if $\sigma \in GF(I)$.
\end{defn}

\begin{prop}[\textrm{\cite[Proposition 2.4]{Stu}}]
The Gr\"obner fan for $I$ is a convex polyhedral fan.
\end{prop}

\begin{rem}
In general, the Gr\"obner fan for $I$ does not consist of strongly convex cones.
\end{rem}

\section{Weight polytope and state polytope}

In the first half of this section, we explain that the normalization of the closure of the torus orbit of a point in a projective space is the toric variety corresponding to some polytope, which is called the weight polytope (Proposition \ref{Oda}).

In the latter half of this section, we consider the Hilbert scheme $\Hilb_h(\P^{r-1})$ being embedded to the projective space $\P^n$ by the Pl\"ucker embedding. Then, the weight polytope corresponding to the normalization of the closure in $\Hilb_h(\P^{r-1})$ of the torus orbit of a Hilbert point $I$ is called the state polytope of $I$. By a theorem of Bayer and Morrison \cite{BM}, the normal fan of the state polytope coincides with the Gr\"obner fan for $I$ (Theorem \ref{gfan}). In conclusion, the normalization of the closure of the torus orbit of a Hilbert point $I$ is the toric variety corresponding to the Gr\"obner fan for $I$. We also prove that the normalization of the closure in $\Hilb_h(\A^{r-1})$ of the torus orbit of $I$ in $\Hilb_h(\A^{r-1})$ is the toric variety corresponding to the Gr\"obner fan for the dehomogenization of $I$ (Corollary \ref{lgfan}).\\

Let $V$ be an $n$ dimensional vector space over an algebraically closed field $K$ and $T$ an algebraic torus of dimension $r$.
Let $\rho:T \rightarrow GL(V)$ be a rational linear representation of $T$ such that $\rho(T)$ contains all scalar multiplications of $GL(V)$. Let the point of $\P(V)$ corresponding to $\mathbf{v} \in V\backslash \{ 0 \}$ will be denoted by $v$. Here $GL(V)$ (and hence $T$ via $\rho$) acts on $\P(V)$ in a natural way. For each $v\in \P(V)$, let $\overline {T\cdot v}$ be the closure in $\P(V)$ of the torus orbit of $v$. Then $\overline {T\cdot v}$ has the open dense orbit which is isomorphic to $T/\textrm{Stab}(v)$.

In general, $\overline{T \cdot v}$ is not a normal variety. But its normalization is the toric variety which contains $T/\mathrm{Stab}(v)$ as an open dense orbit (cf. \cite{Oda}). The corresponding fan of this toric variety is described by the weight polytope explained below.\\

Let $M$ be the character group of $T$ and $N$ the dual lattice of $M$. It is known that $V$ is the direct sum of its weight subspaces:
$$V=\bigoplus_{w \in M} V_{w}, \ \mathrm{where} \ V_w = \{\textbf{v}\in V|\rho(t)\cdot\textbf{v}=w(t)\cdot \textbf{v}\;\mathrm{for\ all}\; t\in T\}.$$

\begin{defn}
Let $\textbf{v}=\sum_{w \in M}\textbf{v}_w$, where $\textbf{v}_w \in V_w$. The convex hull $\mathrm{Wt}(\mathbf{v}) \subset M_\R=M\bigotimes_\Z \R$ of the set $\{w\in M|\textbf{v}_w\neq 0\}$ will be called the \textit{weight polytope} of $\mathbf{v}$.

Let $v$ be a point of $\P(V)$. We denote by $\mathrm{Wt}(v)$ the weight polytope $\mathrm{Wt}(\mathbf{v})$ of a lift $\mathbf{v}$ of $v$. This definition is well-defined.
\end{defn}

When $T \rightarrow T \cdot v$ is injective, the corresponding fan to the normalization of $\overline {T\cdot v}$ is nothing but the normal fan of $\mathrm{Wt}(v)$ (see \cite[Chapter 2.4]{Oda}).

In general, we take a sublattice of $M$ and construct the polytope corresponding to the toric variety as follows.\\

Let $M'$ be the character group of $T/\mathrm{Stab}(v)$ and $N'$ the dual lattice of $M'$. $M'$ is identified with the sublattice of $M$ which is generated by $w_1-w_2$ for any $w_1,w_2 \in \{w\in M|\mathbf{v}_w\neq 0\}$, that is,
$$M'=\left< \{w_1-w_2 \in M | \mathbf{v}_{w_1}, \mathbf{v}_{w_2} \neq 0 \} \right>.$$

The next proposition follows from \cite[Theorem 2.22]{Oda}.

\begin{prop}[cf. \textrm{\cite[Theorem 2.22]{Oda}}]\label{Oda}
The normalization of $\overline {T \cdot v}$ is isomorphic to the toric variety corresponding to the normal fan in $N_\R'$ of the polytope $\mathrm{Wt}(v)-w \subset M_\R'$ for any $w \in \mathrm{Wt}(v) \cap M$.
\end{prop}

\begin{pf}
Let $X$ be the toric variety corresponding to the polytope $\mathrm{Wt}(v)-w$ and let the set of weights $\{w_1,\ldots,w_s\}=\{w\in M|\textbf{v}_w\neq 0\}$. Take a basis of $V_w$ for each $w \in M$ which contains $\mathbf{v_w}$ if $\mathbf{v_w} \neq 0$. Gathering these bases, we get a base of $V$, and we may assume that the homogeneous coordinate of $v$ with respect to this base is given by $(1: \ldots :1:0: \ldots :0)$. We have the morphism $\phi:X \rightarrow \overline{T \cdot v} \subset \P(V)$ as follows:
$$\phi:X \rightarrow \overline{T \cdot v} \subset \P(V) \ \mathrm{;}\  x \mapsto (\chi_{w_1}(x):\ldots:\chi_{w_s}(x):0:\ldots:0),$$
where $\chi_w$ is the character corresponding to $w$. Since $X$ is normal, the morphism $\phi$ factors through the normalization of $\overline{T \cdot v}$. The open torus  orbit of $X$ is isomorphic to the open torus orbit of the normalization of $\overline{T\cdot v}$ by the morphism $\phi$. Moreover, the torus actions on $X$ and on the normalization of $\overline{T \cdot v}$ are compatible with the morphism $\phi$. Then the proposition follows from \cite[Theorem 1.5]{Oda}.
\end{pf}

In the rest of this section, we consider the situation where the torus $T=(K^*)^r$ acts naturally on the vector space $V=N\bigotimes_\Z K \cong K^r$. Then we have a natural embedding $T \subset V$.
Let $S$ be the symmetric algebra of the dual space of $V$, and put $\P=\mathrm{Proj}(S)$. Then $T$ acts on $S=K[x_1, \ldots ,x_r]$ in the following way:
$$T \times S \rightarrow S \ \mathrm{;}\ (t=(t_1,\ldots,t_r),f(X)) \mapsto t \cdot f(X) := f(\frac{x_1}{t_1}, \ldots, \frac{x_r}{t_r}).$$
Further $T$ acts on $\P$, too.

Let $h(x)$ be a polynomial and $\Hilb_h(\P)$ the Hilbert scheme corresponding to $h$. We identify $I \in \Hilb_h(\P)$ with the corresponding homogeneous ideal of $S$. Then $d$-graded part $I_d$ of $I$ is a subspace of $S_d$ of codimension $h(d)$ for any large integer $d \gg 0$.
Let $G(l-h(d),S_d)$ be the Grassmannian variety of subspaces of $S_d$ of codimension $h(d)$, where $l$ is 
$
\dim_KS_d =
\left(
\begin{array}{c}
r+d-1 \\
d \\
\end{array}
\right).
$
It is known that the mapping $\Hilb_h(\P) \rightarrow G(l-h(d),S_d)\ ;\ I \mapsto I_d$ is a closed embedding for any large integer $d \gg 0$.
By the Pl\"ucker embedding, we have $\Hilb_h(\P) \subset G(l-h(d),S_d) \subset \P(\wedge^{l-h(d)}S_d) \cong \P^n$, where $n=\dim_K(\wedge^{l-h(d)}S_d)-1$, and we have a natural action of $T$ on $S_d$, and hence a natural action on $G(l-h(d),S_d)$, and also on $\P^n$.

Then we can apply the facts in the first half of this section to the torus orbit of a Hilbert point $I$.

\begin{defn}\label{state}
Let $I_{d}$ be the $d$-graded component of $I \in \Hilb_h(\P)$. The weight polytope of the Hilbert point $I_d \in \P^n$ is called the \textit{state polytope} $\mathrm{St}_d(I)$ of $I$ in degree $d$.
\end{defn}

The following is deduced from a theorem of Bayer and Morrison \cite{BM} (cf. Sturmfels \cite{Stu2}). For reader's convenience, we give here a slightly different proof.

\begin{thm}[\textrm{\cite[Theorem 2.1]{Stu2}}]\label{gfan}
Let $w \in \mathrm{St}_d(I) \cap M$. Then the normal fan in $N'_\R$ of the polytope $\mathrm{St}_d(I)-w \subset M'_\R$ coincides with the set of the images of the Gr\"obner cones for $I$ by the natural projection $N_\R \rightarrow N_\R'$.
\end{thm}

\begin{pf}
It is enough to prove that the pull back of the normal fan of $\mathrm{St}_d(I)-w$ coincides with the Gr\"obner fan for $I$. Therefore, we prove that the normal fan in $N_\R$ of the state polytope coincides with the Gr\"obner fan for $I$.
First note that each Gr\"obner cone contains the vector $(1,\ldots ,1)$, and the kernel of the projection $N_\R \rightarrow N_\R'$ also contains the vector $(1,\ldots ,1)$. Then we only have to prove that the intersection of the normal fan in $N_\R$ of the state polytope and $\R_+^r \subset N_\R$ coincides with the intersection of the Gr\"obner fan for $I$ and $\R_+^r$.

First, we prove that the latter is a subdivision of the former. We recall that the coordinate of $\P^n$ is given by the Pl\"ucker coordinate. Let $f_1, \ldots , f_k$ be a $K$-basis of $I_d \subset K[X]_d$, where $k=l-h(d)$. Then $f_1 \wedge \cdots \wedge f_k$ defines the Pl\"ucker coordinate of $I_d$. Let $\mathcal{S}$ be the set of weights defined by its all non-zero components. Then, by definition, the state polytope of $I$ in degree $d$ is the convex hull of $\mathcal{S}$.
 Let $\sigma$ be a maximal cone of the Gr\"obner fan for $I$ and $w \in \sigma \cap \R_+^r$ be a weight vector contained in the relative interior of $\sigma$. Let a set of monomials $X^{a_1}, \ldots , X^{a_k}$ be the basis of $\ini_w(I)_d$. Then the weight $v \in M_\R$ corresponding to $X^{a_1} \wedge \cdots \wedge X^{a_k}$ is in $\mathcal{S}$. For $v'( \neq v) \in \mathcal{S}$, $v'$ is defined by wedge product of $X^{a_i'}$, where $X^{a_i'}$ is a term of $ f_i$, and at least one term does not coincide with $X^{a_i}$. Hence we have $\left< w,v \right> > \left< w,v' \right>$ for any $v' \neq v$. This means that $v$ is a vertex of the state polytope and $w$ is contained in the normal cone of $v$. Hence $\sigma$ is contained in the normal cone of $v$. Therefore, the intersection of the Gr\"obner fan for $I$ and $\R_+^r$ is a subdivision of the intersection of the normal fan in $N_\R$ of the state polytope and $\R_+^r$.

Let $\sigma' \neq \sigma$ be a maximal cone of the Gr\"obner fan for $I$ and $w' \in \sigma' \cap \R_+^r$ a weight vector contained in the relative interior of $\sigma'$. Let a set of monomials $X^{b_1}, \ldots ,X^{b_k}$ be the basis of $\ini_{w'}(I)_d$. Then we have $X^{a_1} \wedge \cdots \wedge X^{a_k} \neq X^{b_1} \wedge \cdots \wedge X^{b_k}$. Hence the intersection of the normal fan in $N_\R$ of the state polytope  and $\R_+^r$ coincides with the intersection of the Gr\"obner fan for $I$ and $\R_+^r$.
\end{pf}

In the rest of this section, we consider the case where the dimension of $\mathrm{Stab}(v)$ is one dimensional. Then $\mathrm{Stab}(v)$ is generated by $(1,\ldots,1)K^*$ and a finite abelian group. Therefore, $N_\R'$ is identified with $N_\R/(1,\ldots,1)$.

\begin{cor}\label{lgfan}
Let $I'$ be an ideal of $K[\frac{x_1}{x_r}, \ldots , \frac{x_{r-1}}{x_r}]$ and $I \subset K[x_1, \ldots . x_r]$ the homogenization of $I'$. We assume that the stabilizer of $I$ is one dimensional. We take the image of the part $e_1, \ldots , e_{r-1}$ of the standard basis of $N_\R$ as the basis of $N_\R'$. Then the Gr\"obner fan for $I'$ consists of the images of the Gr\"obner cones for $I$ which contained in the Gr\"obner region $GR(I')$ for $I'$ in $N_\R'$.
\end{cor}

\begin{pf}
We only have to prove that the image $\overline{\sigma}$ of a maximal Gr\"obner cone $\sigma$ for $I$ coincides with some maximal Gr\"obner cone for $I'$ if and only if the intersection of $\R_+^{r-1}\cong \sum_{i=1}^{r-1}\R_+ e_i \subset N_\R'$ and the relative interior of $\overline{\sigma}$ is not empty.

The necessity is clear from the definition. We prove the sufficiency.

For a polynomial $f \in K[x_1, \ldots , x_r]$, we denote by $f'$ the dehomogenization of $f$. For a polynomial $f \in K[\frac{x_1}{x_r}, \ldots ,\frac{x_{r-1}}{x_r}]$, we denote $f^h$ by the honogenization of $f$. Let $\sigma$ be a maximal Gr\"obner cone for $I$. In the following, we denote by $c^i$ the relative interior of a cone $c$. Let $w_1'$ be a weight vector in  $\overline{\sigma}^i \cap \R_+^{r-1}$. Let $\tau$ be the Gr\"obner cone for $I'$ with respect to $w_1'$. Let $w_2'$ be a weight vector in $\tau^i$. Then we can take the reduced Gr\"obner basis $\mathcal{G}=\{g_1, \ldots , g_k \}$ of $I'$ with respect to $<_{w_1'}$ such that $\ini_{w_1'}(g_j)=\ini_{w_2'}(g_j)$ for any $g_j \in \mathcal{G}$.

Choose $w_1,w_2 \in N_\R$ lifts of $w_1',w_2'$. Let $f$ be any homogeneous polynomial in $I$ and $f' \in I'$ dehomogenization of $f$. Let $\phi : N_\R' \rightarrow N_\R$ be the lattice homomorphism as follows:
$$\phi : N_\R' \rightarrow N_\R\  ;\ \overline{e_i} \mapsto e_i .$$
Then we have 
$$\ini_{w_1}(f)=\ini_{\phi(w_1')}(f)=\ini_{w_1'}(f') \cdot x_r^t,$$
where $t=\deg f - \deg \ini_{w_1'}(f')$.

We have $\ini_{w_1'}(f') \in \left<\ini_{w_1'}(g_1), \ldots ,\ini_{w_1'}(g_k) \right>$, then we have $\ini_{w_1}(f) \in \left<\ini_{w_1}(g_1^h), \ldots ,\ini_{w_1}(g_k^h) \right>$. This implies  $g_1^h, \ldots ,g_k^h$ is the reduced Gr\"obner basis of $I$ with respect to $<_{\phi(w_1')}$.

 In the same way, we have $\ini_{w_2}(f) \in \left<\ini_{w_2}(g_1^h), \ldots ,\ini_{w_2}(g_k^h) \right>$. Hence $w_1$ is $I$-equivalent to $w_2$. Thus we have $\tau \subset \overline{\sigma}$.

On the other hand, let $w_3$ be any weight vector in $\sigma^i$. Then we have $\ini_{w_3}(I)$ is generated by $\ini_{\phi(w_1')}(g_1), \ldots , \ini_{\phi(w_1')}(g_k)$. Let $w_3'$ be the image of $w_3$. Let $f$ be a polynomial in $I'$ and $f^h$ the homogenization of $f$. We have $\ini_{w_3'}(f)=(\ini_{w_3}(f^h))'$ and $\ini_{w_3}(f^h)$ is contained in $\left< \ini_{w_3}(g_1^h), \ldots ,\ini_{w_3}(g_k^h) \right>$. Then we have 
\begin{align*}
\ini_{w_3'}(f) \in \left< (\ini_{w_3}(g_1^h))', \ldots ,(\ini_{w_3}(g_k^h))' \right>&= \left< (\ini_{\phi(w_1')}(g_1^h))', \ldots ,(\ini_{\phi(w_1')}(g_k^h))' \right>\\
&=\left< \ini_{w_1'}(g_1), \ldots , \ini_{w_1'}(g_k) \right>.
\end{align*}
In particular, $\ini_{w_3'}(g_i)$ is contained in $\left< \ini_{w_1'}(g_1), \ldots , \ini_{w_1'}(g_k) \right>$, and $\mathcal{G}$ is the reduced Gr\"obner basis of $I'$ with respect to $<_{w_1'}$. Then we have $\ini_{w_3'}(I')= \ini_{w_1'}(I')$. Hence $w_3'$ is $I'$-equivalent to $w_1'$. This implies $\overline{\sigma} \subset \tau$.
\end{pf}

\begin{cor}
Let $I$ be a homogeneous ideal in $\Hilb_h(\P)$. We assume that the stabilizer $\mathrm{Stab}(I)$ of $I$ is one dimensional and that the intersection of the algebraic torus of $\P$ and the closed set $Z$ defined by $I$ is an open dense subset of $Z$. Let $I_i \subset K[\frac{x_1}{x_i},\ldots,\frac{x_r}{x_i}]$ be the dehomogenization of $I$ and $\Sigma_i$ the Gr\"obner fan for $I_i$ with respect to the basis $\{\overline{e_1} , \ldots , \check{\overline{e_i}} , \ldots , \overline{e_r} \}$ of $N_\R'$. Then the set of the images of the Gr\"obner cones for $I$ coincides with the fan obtained as the union $\cup_i \Sigma_i$ of $\Sigma_i$.
\end{cor}

\begin{pf}
This follows from Corollary \ref{lgfan} immediately. Note that the Gr\"obner region $GR(I_i)$ for $I_i$ contains $\sum_{j \neq i}\R_+ e_j \cong \R_+^{r-1}$.
\end{pf}

\section{Main theorem and its applications}
First, we recall the notation and the results in the previous section.
Let $K$ be an algebraically closed field and $T$ the algebraic torus of dimension $r$ which is identified with the diagonal subgroup of $GL(r,K)$. Let $\P^{r-1}$ be the $(r-1)$-dimensional projective space over $K$ and $S \cong k[x_1, \ldots ,x_r]$ the homogeneous coordinate ring of $\P^{r-1}$, so $(x_1:\ldots :x_r)$ is a homogeneous coordinate of $\P^{r-1}$. Then $T$ acts on $\P^{r-1}$ defined by $(t_1, \ldots ,t_r) \cdot (p_1: \ldots :p_r)=(t_1p_1: \ldots :t_rp_r)$, and $T$ acts on $S$ defined by $(t_1, \ldots ,t_r) \cdot f(x_1, \ldots ,x_r)=f(\frac{x_1}{t_1}, \ldots ,\frac{x_r}{t_r})$. Then $T$ acts on the Hilbert scheme $\Hilb^m(\P^{r-1})$ of $m$ points of $\P^{r-1}$.

Let $M$ be the character group of the torus $T$ and $N$ the dual lattice of $M$. Since $T$ acts on $S$, the Gr\"obner fan for an ideal $I \subset S$ is defined in $N_\R$.
Let $v$ be a point of $\Hilb^m(\P^{r-1})$ and $I$ the homogeneous ideal of $S$ corresponding to $v$. Then, by Proposition \ref{Oda}, the normalization of the closure in $\Hilb^m(\P^{r-1})$ of the torus orbit of $v$ is the toric variety corresponding to the state polytope defined in (\ref{state}) for $d \gg 0$. Note that the state polytope is a translate of a subset of $M'_\R$, where $M'$ is the character group of $T/\mathrm{Stab}(v)$. Let $N'$ be the dual lattice of $M'$. We proved that the normal fan in $N_\R'$ of the state polytope coincides with the set of images of the Gr\"obner cones for $I$ by the natural projection $N_\R \rightarrow N'_\R$ (cf. Theorem \ref{gfan}). Then the set of images of the Gr\"obner cones for $I$ becomes a fan which consists of strongly convex cones.\\

Next, following Nakamura \cite{Nak}, we introduce the $G$-orbit Hilbert schemes $\Hilb^G(\P^{r-1})$. 
Let $G$ be a finite abelian subgroup of $GL(r,K)$ of order prime to the characteristic of $K$. Let $\phi :GL(r,K) \rightarrow PGL(r-1,K)$ be the natural projection. Let $m$ be the order of $G$. We assume that $G \cong \phi(G)$ and that $G$ is a subgroup of the diagonal subgroup of $PGL(r-1,K)$. Then we can identify $T/(1, \ldots ,1)K^*$ with an algebraic torus $T'$ of $\P^{r-1}$ and the actions of $G$ and $T$ are commutative. Since $G$ and $T$ also act on $S$, the groups $G$ and $T$ act on the Hilbert scheme $\Hilb^m(\P^{r-1})$.

Let $\Ch^m(\P^{r-1})$ be the Chow variety of $m$ points in $\P^{r-1}$. We have a natural morphism $\phi:\Hilb^m(\P^{r-1}) \rightarrow \Ch^m(\P^{r-1})$, which is called the Hilbert-Chow morphism. Here $G$ acts on $\Hilb^m(\P^{r-1})$ and on $\Ch^m(\P^{r-1})$, and $\phi$ is $G$-equivariant. Therefore we have a natural morphism between their $G$-fixed point sets.
It is known that the $G$-fixed point set of $\Ch^m(\P^{r-1})$ contains $U_r/G$ as a locally closed subset, where $U_r=\{x_r \neq 0\} \subset \P^{r-1}$, and its closure is an irreducible component of the $G$-fixed point set of $\Ch^m(\P^{r-1})$.

\begin{defn}
The $G$-orbit Hilbert scheme $\Hilb^G(\P^{r-1})$ is defined to be the unique irreducible component of the $G$-fixed point set of $\Hilb^m(\P^{r-1})$ which dominates the closure of $U_r/G$ by the map $\phi$.
\end{defn}

\begin{thm}\label{main}
Let $I$ be the G-invariant ideal of $S$ whose zero set is contained in $T' \subset \P^{r-1}$. Then the toric variety corresponding to the fan consisting of the images in $N_\R'$ of the Gr\"obner cones for $I$ is isomorphic to the normalization of the $G$-orbit Hilbert scheme $\Hilb^G(\P^{r-1})$.
\end{thm}

\begin{pf}
By Proposition \ref{Oda} and Theorem \ref{gfan}, it is enough to prove that $\Hilb^G(\P^{r-1})$ is the closure of the torus orbit of $I$ in $\Hilb^m(\P^{r-1})$. 

We can consider that the torus $T'$ is the diagonal subgroup of $PGL(r-1,K)$. If there exist $g,g' \in G$ and $t \in T/(1,\ldots,1) \cong T'$ such that $gt=g't$, then we have $g=g'$ in $T'$. Here $G$ acts freely on the open torus of $\P^{r-1}$.

We prove that the torus orbit of $I$ coincides with the torus contained in $\Hilb^G(\P^{n-1})$.

Let $I'$ be a point of $T\cdot I$, then there exist $t\in T$ such that $I'=t\cdot I$. Since $G \subset T/(1, \ldots ,1)$, $I'$ is also a $G$-fixed point of $\Hilb^m(\P^{r-1})$.

Conversely, let $J$ be an ideal of $S$ defined by a free $G$-orbit whose zero set is contained in $T'$. Then $J$ determines distinct $m$ points $p_1,\ldots ,p_m$ of $T'$ and $I$ determines distinct $m$ points $q_1,\ldots ,q_m$ of $T'$. Take a $t\in T$ satisfying $t\cdot p_1=q_1$. Then we have
$$\{q_1,\ldots ,q_m\}=\{g\cdot q_1|g\in G\}=\{g\cdot t\cdot p_1|g\in G\}=\{t\cdot g\cdot p_1|g\in G\}= t \cdot \{p_1,\ldots,p_m\}.$$
Therefore $J$ is contained in $T\cdot I$.

Since $T\cdot I$ is an open subset of $\Hilb^G(\P^{r-1})$ and $\Hilb^G(\P^{r-1})$ is irreducible, the closure of $T\cdot I$ coincides with $\Hilb^G(\P^{r-1})$.
\end{pf}

\begin{cor}\label{HilbA}
Let $I$ be the G-invariant ideal of $S$ whose zero set is contained in $T'$ and $I' \subset K[\frac{x_1}{x_r}, \ldots , \frac{x_{r-1}}{x_r}]$ the dehomogenization of $I$. Then the toric variety corresponding to the Gr\"obner fan $GF(I')$ for $I'$ is isomorphic to the normalization of $G$-orbit Hilbert scheme $\Hilb^G(\A^{r-1})$.
\end{cor}

\begin{rem}
In this case, the Gr\"obner fan for $I'$ consists of strongly convex cones.
\end{rem}

\begin{pf}
We denote by $\Hilb^G_{\mathrm{norm}}(\P^{r-1})$ (resp. $\Hilb^G_{\mathrm{norm}}(\A^{r-1}$)) the normalization of $\Hilb^G(\P^{r-1})$ (resp. $\Hilb^G(\A^{r-1})$).
A point $p \in \Hilb^G_{\mathrm{norm}}(\P^{r-1})$ is contained in $\Hilb^G_{\mathrm{norm}}(\A^{r-1})$ if and only if $p$ is contained in the torus orbit corresponding to a Gr\"obner cone $\sigma$ of $GF(I')$. Then this corollary follows from Corollary \ref{lgfan} and Theorem \ref{main} immediately.
\end{pf}

\begin{cor}[Ito \textrm{\cite[Theorem 1.1]{Ito}}]
Let $G$ be a finite small cyclic group in $GL(2,\C)$ and $I$ the $G$-invariant ideal of $\C[x,y]$ whose zero set is contained in $(\C^\times)^2$. Then the toric variety corresponding to the Gr\"obner fan for $I$ is isomorphic to the minimal resolution of $\C^2/G$.
\end{cor}
\begin{pf}
Corollary \ref{HilbA} gives that the toric variety corresponding to the Gr\"obner fan for $I$ is isomorphic to the normalization of $G$-orbit Hilbert scheme $\Hilb^G(\A^2)$. Ishii \cite[theorem 3.1]{Ish} shows that the $G$-orbit Hilbert scheme $\Hilb^G(\C^2)$ is the minimal resolution of $\C^2/G$ when $G$ is a finite small subgroup of $GL(2,\C)$. Hence the proof completes.
\end{pf}

\section{Example}
\begin{eg}
Let $G$ be a cyclic group which is generated by the matrix 

$
\left(
\begin{array}{ccc}
e &0 &0  \\
0 &e^{3} &0  \\
0 &0 &1 
\end{array}
\right),
$\\
where $e$ is a primitive fifth root of the unity.
Let $I$ be the ideal of $K[x,y,z]$ which is generated by $x^3-yz^2,x^{2}y-z^3$. The set of zeros of $I$ is $\{ (1:1:1)$,$(e:e^{3}:1)$,$(e^{2}:e^{1}:1)$,$(e^{3}:e^{4}:1),(e^{4}:e^{2}:1)\}$. Then $I$ is a $G$-invariant ideal whose zero set is contained in $(K^*)^2$.

We denote by $\mathcal{G}_w$ the reduced Gr\"obner basis with respect to $w$ and $c[w]$ the Gr\"obner cone with respect to $w$.
The Gr\"obner fan for $I$ is defined by 11 maximal cones:

\begin{figure}[h]
\begin{center}
\includegraphics[width=15em]{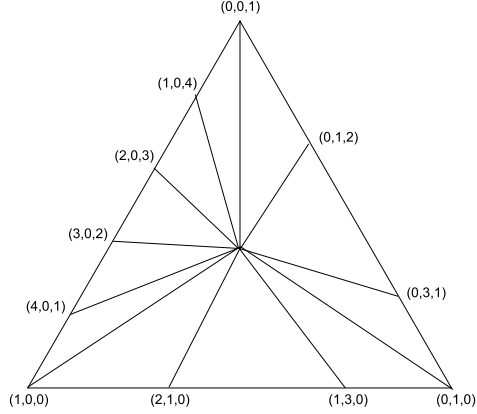}
\end{center}
\caption[defn]{The Gr\"obner fan for $I$}
\label{fig-label-1}
\end{figure}

\begin{align*}
\mathcal{G}_{w_1}&=\{ x^3-yz^2,x^2y-z^3,xz^3-y^2z^2,xy^3z^2-z^6,y^5z^2-z^7 \}\\
\mathcal{G}_{w_2}&=\{ x^3-yz^2,x^2y-z^3,y^2z^2-xz^3 \}\\
\mathcal{G}_{w_3}&=\{ yz^2-x^3,x^2y-z^3,x^5-z^5 \}\\
\mathcal{G}_{w_4}&=\{ yz^2-x^3,x^2y-z^3,z^5-x^5 \}\\
\mathcal{G}_{w_5}&=\{ yz^2-x^3,z^3-x^2y,x^2y^2-x^3z \}\\
\mathcal{G}_{w_6}&=\{ yz^2-x^3,z^3-x^2y,x^3z-x^2y^2 \}\\
\mathcal{G}_{w_7}&=\{ yz^2-x^3,z^3-x^2y,x^3z-x^2y^2,x^2y^3z-x^6,x^7-x^2y^5 \}\\
\mathcal{G}_{w_8}&=\{ yz^2-x^3,z^3-x^2y,x^3z-x^2y^2,x^6-x^2y^3z \}\\
\mathcal{G}_{w_9}&=\{ x^3-yz^2,z^3-x^2y,y^2z^2-xz^3 \}\\
\mathcal{G}_{w_{10}}&=\{ x^3-yz^2,x^2y-z^3,xz^3-x^2y^2,z^6-xy^3z^2 \}\\
\mathcal{G}_{w_{11}}&=\{ x^3-yz^2,x^2y-z^3,xz^3-y^2z^2,xy^3z^2-z^6,z^7-y^4z^2 \}
\end{align*}
\begin{align*}
c[w_1]&=\{(x,y,z) \in \R^3|y-z\geq0,x-2y+z\geq0\}\\
c[w_2]&=\{(x,y,z) \in \R^3|x-2y+z\leq0,3x-y-2z\geq0\}\\
c[w_3]&=\{(x,y,z) \in \R^3|3x-y-2z\leq0,x-z\geq0\}\\
c[w_4]&=\{(x,y,z) \in \R^3|x-z\leq0,2x+y-3z\geq0\}\\
c[w_5]&=\{(x,y,z) \in \R^3|2x+y-3z\leq0,x-2y+z\leq0\}\\
c[w_6]&=\{(x,y,z) \in \R^3|x-2y+z\geq0,x-y\leq0\}\\
c[w_7]&=\{(x,y,z) \in \R^3|x-y\geq0,4x-3y-z\leq0\}\\
c[w_8]&=\{(x,y,z) \in \R^3|4x-3y-z\geq0,3x-y-2y\leq0\}\\
c[w_9]&=\{(x,y,z) \in \R^3|3x-y-2z\geq0,2x+y-3z\leq0\}\\
c[w_{10}]&=\{(x,y,z) \in \R^3|2x+y-3z\geq0,x+3y-4z\leq0\}\\
c[w_{11}]&=\{(x,y,z) \in \R^3|x+3y-4z\geq0,y-z\leq0\}
\end{align*}
Let $N$ be the lattice corresponding to the torus $T$. The stabilizer of $I$ is generated by $(1,1,1)\cdot K^{*}$ and $(e,e^3,1)$. Then the lattice $N'$ corresponding to $T/\mathrm{Stab}(I)$ is identified with $\Z\oplus\Z\oplus(\tfrac{1}{5},\tfrac{3}{5})\Z$ and the lattice homomorphism $N \rightarrow N'$ is described as follows:
$$N\cong\Z^3 \rightarrow N'\cong\Z\oplus\Z\oplus(\tfrac{1}{5},\tfrac{3}{5})\Z\  \mathrm{:}\ (a,b,c) \mapsto (a-c,b-c).$$
Then the fan corresponding to $\Hilb^G(\P^2)$ is generated by the images of these cones.

\begin{figure}[h]
\begin{center}
\includegraphics[width=15em]{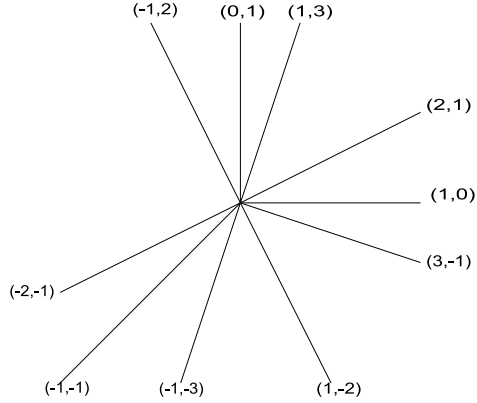}
\end{center}
\caption[defn]{The image of the Gr\"obner fan for $I$}
\label{fig-label-1}
\end{figure}

This fan corresponds to a nonsingular toric variety, and the subfan which is contained in $\R_+^2$ corresponds to the minimal resolution of $\{(x:y:z) \in \P^2 | z \neq 0 \}/G \cong  K^2/G$.
\end{eg}
\begin{eg}
Let $G$ be a cyclic group which is generated by the matrix

$
\left(
\begin{array}{cccc}
e &0 &0 &0  \\
0 &e^2 &0 &0  \\
0 &0 &e^3 &0  \\
0 &0 &0 &1 
\end{array}
\right),
$\\
where $e$ is a primitive fifth root of the unity.
Let $I$ be the ideal of $K[x,y,z,w]$ which is generated by $x^5-w^5,x^{2}-yw,x^3-zw^2$. The set of zeros of $I$ is $\{ (1:1:1:1),(e:e^2:e^3:1),(e^2:e^4:e:1),(e^3:e:e^4:1),(e^4:e^3:e^2:1)\}$. Then $I$ is a $G$-invariant ideal whose zero set is contained in $(K^*)^3$. Here $G$ acts on $\{w \neq 0 \} \cong \A^3$. Let $I' \subset K[\frac{x}{w},\frac{y}{w},\frac{z}{w}] $ be the dehomogenization of $I$. The Gr\"obner fan for $I'$ has 15 edges, 32 facets, and 18 maximal cones. 17 maximal cones are simplicial and nonsingular, but 1 maximal cone is not simplicial. Then $\Hilb^G(\A^3)$ is singular.
\end{eg}

\begin{figure}[h]
\begin{center}
\includegraphics[width=15em]{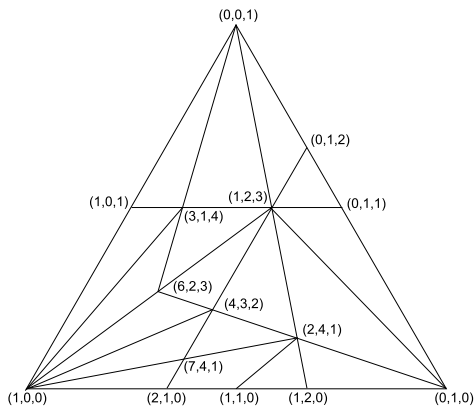}
\end{center}
\caption[defn]{The Gr\"obner fan for $I'$}
\label{fig-label-1}
\end{figure}

Tomohito Morita\\
Department of Mathematics\\
Tokyo Institute of Technology\\
Oh-okayama, Meguro-ku, Tokyo 152-8551\\
Japan\\
morita.t.ae@m.titech.ac.jp

\begin{thebibliography}{99}
\bibitem{BM}
D.\ Bayer and I.\ Morrison, Gr\"obner bases and geometric invariant theory I. Initial ideals and state polytopes, J.\ Symbolic Computation 6 (1988), 209--217.
\bibitem{CMT}
A.\ Craw, D.\ Maclagan, and R.\ R.\ Thomas, Moduli of McKay quiver representation I: The coherent component, Proc.\ London Math.\ Soc.\ (2007), 95 (1):179--198.
\bibitem{CMT2}
A.\ Craw, D.\ Maclagan, and R.\ R.\ Thomas, Moduli of McKay quiver representation II: Gr\"obner basis techniques, to appear in J.\ Algebra.
\bibitem{Ful}
W.\ Fulton, Introduction to Toric Varieties, Princeton University Press.
\bibitem{Jen}
A.\ Jensen, Gfan - a software system for Gr\"obner fans. Available at

http://home.imf.au.dk/ajensen/software/gfan/gfan.html 
\bibitem{KSZ}
M.\ M.\ Kapranov, B.\ Sturmfels and A.\ V.\ Zelevinsky, Chow polytopes and general resultants, Duke Math.\ J.\ Volume 67, Number 1 (1992), 189--218. 
\bibitem{Ish}
A.\ Ishii, On the McKay correspondence for a finite small subgroup of $GL(2,\mathbb{C})$, J.\ Reine Angew.\ Math.\ 549 (2002), 221--233.
\bibitem{Ito}
Y.\ Ito, Minimal resolution via Gr\"obner basis, Algebraic Geometry in East Asia, (IIAS, 2001), World Scientific, (2003), 165--174.
\bibitem{IN}
Y.\ Ito and I.\ Nakamura, Hilbert schemes and simple singularities, In New trends in algebraic geometry (Warwick, 1996), volume 264 of London Math.\ Soc.\ Lecture Note Ser., pages 151--233. Cambridge Univ.\ Press, Cambridge, 1999.
\bibitem{IN2}
Y.\ Ito and I.\ Nakamura, McKay correspondence and Hilbert schemes, Proc.\ Japan Acad.\ 72 (1996), 135--138.
\bibitem{Nak}
I.\ Nakamura, Hilbert schemes of abelian group orbits, J.\ Algebraic Geom., 10 (4) : 757--779, 2001.
\bibitem{Oda}
T.\ Oda, Convex Bodies and Algebraic Geometry, Ergeb.\ Math.\ Grenzeb. (3) 15, Springer-Verlag, Berlin, 1988.
\bibitem{Stu}
B.\ Sturmfels, Gr\"obner Bases and Convex Polytopes, Univ.\ Lect.\ Series, 8, AMS (1995).
\bibitem{Stu2}
B.\ Sturmfels, Gr\"obner bases of toric varieties, Tohoku Math.\ J.\ 43 (1991), 249--261.
\bibitem{Yas}
T.\ Yasuda, Universal flattening of Frobenius, arXiv:0706.2700v3 [math.AG]
\end{thebibliography}
\end{document}